\documentclass[10pt]{amsart}

\usepackage{a4wide}
\usepackage[utf8]{inputenc}
\usepackage{amsmath, amssymb, bbm}
\usepackage{amsthm}
\usepackage{mathtools}
\usepackage{bbm}
\usepackage{blkarray}
\usepackage[matrix,arrow,curve]{xy}
\usepackage{tikz}
\usetikzlibrary{cd}
\usepackage{graphicx}
\usepackage{comment}
\usepackage{color}
\usepackage[normalem]{ulem}
\usepackage{enumerate}
\usepackage{enumitem}
\definecolor{darkblue}{rgb}{0,0,0.6}
\usepackage[ocgcolorlinks,colorlinks=true, citecolor=darkblue, filecolor=darkblue, linkcolor=darkblue, urlcolor=darkblue]{hyperref}
\usepackage[capitalize,noabbrev]{cleveref}
\usepackage[mathscr]{euscript}
\usepackage{stmaryrd}

\usepackage{sansmath}

\newcommand{\cocolon}{\nobreak \mskip6mu plus1mu \mathpunct{}\nonscript\mkern-\thinmuskip {:}\mskip2mu \relax}

\newcommand*{\coloneqq}{\mathrel{\vcenter{\baselineskip0.5ex \lineskiplimit0pt
                     \hbox{\scriptsize.}\hbox{\scriptsize.}}}%
                     =}
                     
\usepackage{pict2e,picture}
\makeatletter
\DeclareRobustCommand{\bbDelta}{{\mathpalette\bb@Delta\relax}}
\newcommand{\bb@Delta}[2]{%
  \begingroup
  \sbox\z@{$\m@th#1\Delta$}%
  \dimendef\Dht=6 \dimendef\Dwd=8
  \setlength{\Dwd}{\wd\z@}%
  \setlength{\Dht}{\ht\z@}%
  \begin{picture}(\Dwd,\Dht)
  \put(0,0){$\m@th#1\Delta$}
  \put(.52\Dwd,.77\Dht){\line(-13,-26){.35\Dht}}
  \end{picture}%
  \endgroup
}
\makeatother

\DeclareFontFamily{T1}{cbgreek}{}
\DeclareFontShape{T1}{cbgreek}{m}{n}{<-6>  grmn0500 <6-7> grmn0600 <7-8> grmn0700 <8-9> grmn0800 <9-10> grmn0900 <10-12> grmn1000 <12-17> grmn1200 <17-> grmn1728}{}
\DeclareSymbolFont{quadratics}{T1}{cbgreek}{m}{n}
\DeclareMathSymbol{\qoppa}{\mathord}{quadratics}{19}
\DeclareMathSymbol{\Qoppa}{\mathord}{quadratics}{21}

\newtheoremstyle{introthms}
	{}{}{\itshape}{}{\bfseries }{}{ }
	{\thmname{#1} \thmnumber{#2}. \thmnote{\bfseries{(#3)}}}

\newtheoremstyle{thms}
	{}{}{\itshape}{}{\bfseries }{}{ }
	{\thmname{#1}\thmnumber{#2}. \thmnote{\bfseries{[#3]}}}
\newtheoremstyle{thms2}
	{}{}{\itshape}{}{\bfseries }{}{ }
	{\thmname{#1}\thmnumber{#2}. \thmnote{\bfseries{[#3]}}}

\newtheoremstyle{name}
	{}{}{\itshape}{}{\bfseries }{}{ }
	{\thmname{#1}\thmnumber{#2}\thmnote{\bfseries{[#3]}}}
\newtheoremstyle{defs}
	{}{12pt}{\normalfont}{}{\bfseries }{}{ }
	{\thmname{#1} \thmnumber{#2}. \thmnote{\bfseries{(#3)}}}
\newtheoremstyle{defs2}
	{}{12pt}{\normalfont}{}{\bfseries }{}{ }
	{\thmname{#1}\thmnumber{#2}. \thmnote{\bfseries{(#3)}}}
\newtheoremstyle{rmk}
	{}{}{\normalfont}{}{\itshape }{}{ }{\thmname{#1}. \thmnote{#3}}
\newtheoremstyle{claim}
	{}{}{\normalfont}{}{\itshape}{}{ }{\thmname{#1} \thmnumber{#2}. \thmnote{#3}}

\theoremstyle{introthms}

\swapnumbers
	
\theoremstyle{thms}
\newtheorem{proposition}{Proposition}[section]
\newtheorem{lemma}[proposition]{Lemma}%[section]

\newtheorem*{uthma}{Theorem A}
\newtheorem*{uthmb}{Theorem B}

\theoremstyle{defs}
\newtheorem{definition}[proposition]{Definition}%[section]

\newtheorem{example}[proposition]{Example}
\newtheorem{remark}[proposition]{Remark}

\theoremstyle{defs2}

\theoremstyle{rmk}

\theoremstyle{claim}

\newcommand{\op}{\mathrm{op}}              % opposite category

\newcommand{\fib}{\mathrm{fib}}            % homotopy fibre
\newcommand{\cof}{\mathrm{cof}}            % homotopy quotient
\renewcommand{\hom}{\mathrm{hom}}          % spectral Hom functor
\newcommand{\colim}{\operatorname{colim}}

\newcommand{\Cat}{\mathrm{Cat}}      
      
\newcommand{\Fun}{\mathrm{Fun}}            % functor cat
\newcommand{\Ek}{{\mathbb E_k}}
\newcommand{\Alg}{{\mathrm{Alg}}}

\newcommand{\Dperf}{\mathcal D^{\mathrm{p}}}

\newcommand{\Mod}{{\mathrm{Mod}}}
\newcommand{\Sp}{\mathrm{Sp}}
\newcommand{\GEM}{\mathrm{H}}

\newcommand{\C}{\mathcal C}                % stable infty cat
\newcommand{\D}{\mathcal D}                % stable infty cat
\newcommand{\id}{\mathrm{id}}              % identity

\newcommand{\Eone}{\mathbb E_1}

\newcommand{\Ar}{\operatorname{Ar}}

\makeatletter
\providecommand{\leftsquigarrow}{%
  \mathrel{\mathpalette\reflect@squig\relax}%
}
\newcommand{\reflect@squig}[2]{%
  \reflectbox{$\m@th#1\rightsquigarrow$}%
}
\makeatother

\newcommand{\Einf}{\mathbb E_\infty}

\title{A note on higher almost ring theory}

\author[Fabian Hebestreit]{Fabian Hebestreit}
\author[Peter Scholze]{Peter Scholze}

\makeatletter
\let\@wraptoccontribs\wraptoccontribs
\makeatother

\date{}
\begin{document}
\setcounter{tocdepth}{1}

\begin{abstract}
We explain a derived version of the basic construction of localisations of module categories by means of idempotent ideals, which lie at the heart of Faltings' almost ring theory. We use it to provide an example of a commutative $\mathbb F_p$-algebra whose Frobenius endomorphism does not induce an isomorphism on its smashing spectrum.
\end{abstract}

\maketitle
\tableofcontents

%%%%%%%%%%%%%%%%%%%%%%%%%%%%%%%%%%%%%%%%%%%%%%%%%%%%%%%%%%%%%%%%
%%%%%%%%%%%%%%%%%%%%%%%%%%%%%%%%%%%%%%%%%%%%%%%%%%%%%%%%%%%%%%%%
%%%%%%%%%%%%%%%%%%%%%%%%%%%%%%%%%%%%%%%%%%%%%%%%%%%%%%%%%%%%%%%%

\section{Introduction}

Almost ring theory was introduced by Faltings in \cite{Fal1, Fal2}, as a way capturing and propagating vanishing phenomena in Galois cohomology, building on initial work of Tate in \cite{Tate}. The basic set-up of the theory was then reworked by Gabber and Ramero in \cite{GR, GR2} to simply rely on a commutative ring $R$ and ideal $I \subseteq R$ satisfying two assumptions:
\begin{enumerate}
\item $I$ is idempotent (that is $I = I^2$), and
\item $I$ is flat as an $R$-module.
\end{enumerate}
The most prominent example of this situation is given by the ideal of topologically nilpotent elements $I$ inside the ring of power bounded elements $R$ of a perfectoid field. %Gabber and Ramero in fact showed that, without much loss, the second condition can be generalised further by only requiring $I \otimes_R I$ to be flat over $R$, and eventually that one can do largely without it.

For an idempotent ideal $I$, one says that a morphism of $R$-modules is an \emph{$I$-almost isomorphism} if its kernel and cokernel $I$-almost vanish, that is they are annihilated by all elements of $I$. If $I$ is flat, the localisation $\mathrm{aMod}_I(R)$ of $\Mod(R)$ at these maps, the category of \emph{$I$-almost $R$-modules}, retains many good homological properties: For example the tensor product of $R$-modules descends to it and its derived category can be described as the localisation of the derived category of $R$-modules localised at those maps inducing $I$-almost isomorphisms on homology groups. In fact, these desirable properties are all direct consequences of the fact that the multiplication 
\[R/I \otimes^\mathbb L_R R/I \longrightarrow R/I\] 
is an equivalence, if $I$ satisfies the properties listed above, making $R \rightarrow R/I$ into what is sometimes called a derived localisation.

The purpose of the present note is to explain that by passing to derived categories directly one can do away with the flatness hypothesis  while still retaining this simple explanation for the good properties of the derived category of $I$-almost modules.
Namely, Gabber and Ramero already showed with some effort that the subcategory of the derived category $\D(R)$ spanned by those complexes whose homology $I$-almost vanishes admits both adjoints the moment $I$ is idempotent. From a more modern point of view this also follows from the reflection principle of \cite{reflection}, since this subcategory is easily checked closed under limits and colimits. The Schwede-Shipley theorem now implies that this subcategory is still the derived category of some $R$-algebra in $\D(R)$, that is a derived localisation of $R$. Gabber and Ramero already observed that when $I \otimes_R I$ (but not necessarily $I$) is flat over $R$, this is given by the commutative differential graded algebra $R/^{\mathbb L} (I \otimes_R I)$ with % given by
%\[[R/^{\mathbb L} (I \otimes_R I)]_i = \begin{cases} R & i = 0 \\ I \otimes_R I  & i=1 \\0 & \text{else}\end{cases}\]

\[\mathrm H_i(R/^{\mathbb L} (I \otimes_R I)) = \begin{cases} R/I & i = 0 \\ \mathrm{ker}(I \otimes_R I \rightarrow I) & i = 1 \\ 0 & \text{else.} \end{cases}\]
Our main results identify it in general, even when the base ring is not assumed static. Furthermore, our arguments do not rely on the reflection principle or an a priori analysis of the category of complexes with $I$-almost vanishing cohomology. Namely, we directly show:

\begin{uthma}
\label{theorem:almoststuff}%
Let $R$ be an animated commutative ring
% connective $\Ek$-ring $A$ with $1 \leq k \leq \infty$, respectively. 
and consider the full subcategory $\mathrm{LQ}_R$ of $R/\mathrm{AnCRing}$ spanned by the maps $\varphi \colon R \rightarrow S$ for which
\begin{enumerate}
\item the multiplication $S \otimes^\mathbb L_R S \rightarrow S$ is an equivalence, i.e.\ $\varphi$ is a derived localisation, and
%\item $B$ is connective, and
\item $\pi_0(\varphi) \colon \pi_0R \rightarrow \pi_0S$ is surjective.
\end{enumerate}
Then the functor
\[
\mathrm{LQ}_R \longrightarrow \{I \subseteq \pi_0R \mid I^2 = I\}, \quad \varphi \longmapsto \ker(\pi_0\varphi)
\]
is an equivalence of categories, where we regard the target as a poset via the inclusion ordering. The inverse image of some $I \subseteq \pi_0(R)$ is given by the limit of the Amitsur complex for the map $R \rightarrow \pi_0(R)/I$. 
\end{uthma}

In fact, there is an entirely analogous result for $\mathbb E_k$-rings, which in particular gives a sensible way of forming quotients by idempotent ideals; recall that directly forming quotients even by single elements is neither easy nor generally possible in this realm, see e.g.\ \cite{Burk}:

\begin{uthmb}
Let $A$ be a connective $\Ek$-ring with $1 \leq k \leq \infty$, respectively. Consider again the full subcategory $\mathrm{LQ}_A$ of $A/\mathrm{Alg}_{\mathbb{E}_k}(\Sp)$ spanned by the maps $\varphi \colon A \rightarrow B$ for which
\begin{enumerate}
\item the multiplication $B \otimes_A B \rightarrow B$ is an equivalence, i.e.\ $\varphi$ is a localisation,
\item $B$ is connective, and
\item $\pi_0(\varphi) \colon \pi_0A \rightarrow \pi_0B$ is surjective.
\end{enumerate}
Then the functor
\[
\mathrm{LQ}_A \longrightarrow \{I \subseteq \pi_0A \mid I^2 = I\}, \quad \varphi \longmapsto \ker(\pi_0\varphi)
\]
is an equivalence of categories, where we again regard the target as a poset via the inclusion ordering. The inverse image of some $I \subseteq \pi_0(A)$ can be described more directly as $A / I^\infty$, where 
\[
I^\infty = \lim_{n \in \mathbb N^\op} \widetilde I^{\otimes_A n}
\]
with $\widetilde I \rightarrow A$ the fibre of the canonical map $A \rightarrow \GEM(\pi_0(A)/I)$. Furthermore, this inverse system stabilises on $\pi_i$ for $n > i+1$.
\end{uthmb}

The inverse limit in the statement is formed over the system of maps $\widetilde I^{\otimes_A n+1} \rightarrow \widetilde I^{\otimes_A n}$ induced by the inclusion $\widetilde I \rightarrow A$ in any of the $n+1$ spots; the induced maps agree on account of $\widetilde I$ being a Smith-ideal in $A$, a concept we briefly recall in Section 2 on accout of a lack of a simple general reference.

To connect to the discussion of almost modules before, recall first that for $R$ a commutative animated (e.g.\ static) ring, the derived category of  $R$ depends only on its underlying $\mathbb E_1$-ring $\GEM R$, that is we have $\mathcal D(R) \simeq \Mod(\GEM R)$, and note that the animated commutative ring $S$ corresponding to some idempotent $I \subseteq \pi_0(R)$ necessarly satisfies $\GEM S \simeq (\GEM R)/I^\infty$ by the uniqueness assertions of the theorems above. We shall therefore denote this animated commutative ring by $R/I^\infty$ as well and restrict the discussion to the case of $\mathbb E_1$-rings from here on. Recall then that a localisation of $\mathbb E_1$-rings $\varphi \colon A \rightarrow B$, gives rise to a stable recollement
\[
\begin{tikzcd}
[column sep=7ex]
\Mod(B) \ar[r] 
& \Mod(A) 
\ar[r] 
\ar[l,bend left=30,shift left=1.5ex,start anchor=west,end anchor=east,"{\hom_A(B,-)}"] 
\ar[l,bend right=30,shift right=1.5ex,start anchor=west,end anchor=east,"{B \otimes_A -}"'] 
%\ar[l,phantom,shift left=1.2ex,start anchor=west,end anchor=east,"\perp"] 
%\ar[l,phantom,shift right=1.2ex,start anchor=west,end anchor=east,"\perp"]
& \Mod(A)[\varphi\text{-eq's}^{-1}]
\ar[l,bend left=30,shift left=1.5ex,start anchor=west,end anchor=east,"{\hom_A(\fib(\varphi),-)}"]
\ar[l,bend right=30,shift right=1.5ex,start anchor=west,end anchor=east,"{\fib(\varphi) \otimes_A -}"']
%\ar[l,phantom,shift left=1.2ex,start anchor=west,end anchor=east,"\perp"] 
%\ar[l,phantom,shift right=1.2ex,start anchor=west,end anchor=east,"\perp"]
\end{tikzcd}
\]
where the $\varphi$-equivalences are those maps of $A$-modules whose fibres lie in image of the restriction functor $\Mod(B) \rightarrow \Mod(A)$, see e.g.\ \cite[Appendix A.4]{92}; the diagram above indicates four adjunctions with left adjoints on top, arranged into three horizontal Verdier sequences. 

In the case at hand, $B=A/I^\infty$ with $A$ connective, we now indeed show that the image of the fully faithful restriction $\Mod(A/I^\infty) \rightarrow \Mod(A)$ consists exactly of the $I$-almost vanishing $A$-modules $M$, i.e.\ those with $I \cdot \pi_n(M) = 0$ for all $n \in \mathbb Z$, and consequently that a map is a $\varphi$-equivalence if and only if it induces an $I$-almost isomorphism on all homotopy groups. If $A$ is not assumed connective the same then holds for 
\[A/I^\infty \coloneq A \otimes_{\tau_{\geq 0}A} (\tau_{\geq 0} A)/I^\infty.\]
Keeping the notation $I^\infty$ also for the fibre of $A \rightarrow A/I^\infty$ in the general case, the recollement takes the form

%B \tensor_A cofib(M --> N)

%The kernel of the map $\phi_! \colon \Mod(A) \rightarrow \Mod(A/I^\infty)$ can therefore also be described as the Verdier quotient of $\Mod(A)$ by its full subcategory $\mathrm{Ann}_I(A)$ of those modules, whose homotopy is annihilated by $I$. This quotient is the $\infty$-category of \emph{almost $A$-modules} with respect to $I$ and we denote it by $\mathrm{aMod}_I(A)$. Using this description, the split Verdier sequence associated to $A \rightarrow A/I^\infty$ takes the form
\[
\begin{tikzcd}
[column sep=7ex]
\mathrm{Mod}(A/I^\infty) \ar[r] 
& \Mod(A) 
\ar[r] 
\ar[l,bend left=30,shift left=1.5ex,start anchor=west,end anchor=east,"{\hom_A(A/I^\infty,-)}"] 
\ar[l,bend right=30,shift right=1.5ex,start anchor=west,end anchor=east,"{A/I^\infty \otimes_A -}"'] 
%\ar[l,phantom,shift left=1.2ex,start anchor=west,end anchor=east,"\perp"] 
%\ar[l,phantom,shift right=1.2ex,start anchor=west,end anchor=east,"\perp"]
& \mathrm{aMod}_I(A)
\ar[l,bend left=30,shift left=1.5ex,start anchor=west,end anchor=east,"{\hom_A(I^\infty,-)}"]
\ar[l,bend right=30,shift right=1.5ex,start anchor=west,end anchor=east,"{I^\infty \otimes_A -}"']
%\ar[l,phantom,shift left=1.2ex,start anchor=west,end anchor=east,"\perp"] 
%\ar[l,phantom,shift right=1.2ex,start anchor=west,end anchor=east,"\perp"]
\end{tikzcd}
\]
exhibiting the $I$-almost $A$-modules as a split Verdier quotient of $\Mod(A)$, as desired. Furthermore, if $A$ is an $\mathbb E_k$-algebra, the entire recollement is suitably $\mathbb E_{k-1}$-multiplicative on general grounds.

We end this introduction with three examples, which illustrate the extended range of applicability offered by the removal of the flatness assumption: 
\begin{enumerate}
\item On the one hand, consider for $K$ a field the commutative ring
\[R_n= K[T_1^{1/2^\infty}, \dots, T_n^{1/2^\infty}] \coloneqq K[T_{i,j}, i,j \in \mathbb N, j\leq n]/(T_{i+1,j}^2 - T_{i,j}, i,j \in \mathbb N, j \leq n),\]
together with the ideal 
\[I_n = (T_1^{1/2^\infty}, \dots, T_n^{1/2^\infty}) \coloneqq (T_{i,j}, i,j \in \mathbb N, j \leq n).\]
which is evidently idempotent. It is flat only for $n=1$ but nevertheless $I_n \otimes^\mathbb L_{R_n} I_n = I_n$ so that 
\[R_n/I_n^\infty = R_n/I_n = K\]
is still static. 
\item On the other hand, in $\overline R_n = R_n/(T_1,\dots,T_n)$ the ideal $\overline I_n = I_n/(T_1, \dots, T_n)$ is no longer flat even for $n=1$ and 
\[\pi_*\overline R_n/\overline I_n^\infty = \Lambda_K[T_1, \dots T_n]\]
is the exterior algebra on $n$ generators in degree $1$: The module $\overline I_1 \otimes_{\overline R_1} \overline I_1$ is still flat over $\overline R_1$, giving this calculation for $n=1$ and the general case then follows by multiplicativity. In this example, a well-behaved almost theory thus still exists at the level of modules for $n=1$. This fails for $n>1$, but nevertheless, the derived theory is largely unaffected.
\item On the third hand, our results can be used to analyse wilder examples: For $S$ the set of finite strings of $0$'s and $1$'s take
\[
R=\mathbb F_p[T_s\mid s \in S]/(T_s - T_{s \ast 0} \cdot T_{s \ast 1},T_s^p\mid s \in S)
\]
with $I$ generated by all the variables. Then the maps $R \rightarrow R/I = \mathbb F_p$ and $R \rightarrow R$ have the same extension of scalars under the Frobenius of $R$, and thus so do $R \rightarrow R/I^\infty$ and $R \rightarrow R$ by Theorem A. 
This shows that contrary to the case of the ordinary spectrum, the Frobenius map need not induce the identity or even an isomorphism on the smashing spectrum (which is nothing but the posed of derived localisations) of $R$.
\end{enumerate}

\subsection*{Acknowledgements}

We thank the anonymous referee of \cite{92}, Lars Hesselholt, and Akhil Mathew for making us aware of the connection between almost ring theory and the localisation theorem in algebraic $K$-theory. The clarification of this connection led to the present note, which originally appeared as part of the appendix of \cite{92}, but eventually outgrew it. We further thank Dustin Clausen, Manuel Hoff, Maxime Ramzi and Christoph Winges for several helpful discussions.% about localisations of $\mathbb E_1$-rings, and Manuel Hoff 

The authors were supported by the German Research Foundation (DFG) during the writing of this note, FH through the collaborative research centre “Integral structures in Geometry and Representation theory” (grant no. TRR 358–491392403) at the University of Bielefeld, and PS by a Leibniz Prize, and through the Hausdorff Centre of Mathematics (grant no. EXC-2047/1 – 390685813) at the MPIM Bonn.

\section{Smith ideals}

In the proof of the main results, we will have to pass back and forth between ring maps and their fibres, regarded as ideals. We provide here the basic observations concerning Smith-ideals that facilitate this. The statements are certainly well-known, though we are unaware of a general reference, but see e.g.\ \cite{Hovey, WhiteYau} for a treatment in model categorical language 

Consider then an $\Ek$-monoidal category $(\C,\otimes,\mathbb I)$, $k \geq 1$, and give $\Ar(\C)$ the induced $\Ek$-monoidal Day convolution structure with respect to taking minima on $[1]$ (which is symmetric monoidal). If $\C$ has finite colimits which are respected by the monoidal structure this is again an $\Ek$-monoidal category. Explicitly, we then have
\[(c \rightarrow d) \otimes^\mathrm{Day} (c' \rightarrow d') \simeq (c \otimes d' +_{c \otimes c'} d \otimes c' \rightarrow d \otimes d'),\]
and the unit is $0 \rightarrow \mathbb I$. Even without the assumption, the evaluation functor $t \colon \Ar(\C) \rightarrow \C$ is a map of operads over $\mathbb E_k$, which allows us to set:

\begin{definition}
An \emph{$\Ek$-Smith-ideal} in an $\Ek$-algebra $A$ in $(\C,\otimes,\mathbb I)$ is an $\Ek$-algebra $J \rightarrow A$ in $\Ar(\C)$ lifting the $\Ek$-structure on $A$. Let us denote the category of Smith-ideals in $A$, that is the fibre of $t \colon \mathrm{Alg}_{\mathbb E_k}(\Ar(\C)) \rightarrow \mathrm{Alg}_{\mathbb E_k}(\C)$ over $A$, by $\mathrm{SmId}_{\Ek}(A)$. 
\end{definition}

If $\C$ has finite colimits that are preserved by the monoidal structure, the (source part of the) multiplication map of a Smith-ideal takes the form
\[J \otimes A +_{J \otimes J} A \otimes J \longrightarrow J.\]
It equips $J$ with a non-unital $\mathbb E_k$-structure and an extension of its multiplication to an $A$-$A$-bimodule structure, such that the map $J \rightarrow A$ extends to both a non-unital $\Ek$-algebra and a bimodule map by the following:

\begin{remark} Let $(\D,\otimes)$ be an $\Ek$-monoidal category. There is a general adjunction
\[- \times_{\Ek} \mathcal \D^\otimes \colon \mathrm{Opd}/\Ek \leftrightarrows \mathrm{Opd}/\Ek \cocolon \Fun(\mathcal D^\otimes, -)^{\Ek\text{-}\mathrm{Day}}\]
on the category of operads over $\Ek$. Applying this to $(\mathcal D,\otimes)$ the $\Ek$-monoidal category induced by $([1],\min)$, whose underlying operad is just $\Ek \times [1]^\mathrm{min}$, we see that an $\Ek$-Smith ideal is equivalent to an algebra over the operad $\mathrm{SmId}_{\Ek} \coloneq \Ek \times [1]^\mathrm{min}$. In particular, an $\Ek$-Smith ideal is simply a lax $\Ek$-monoidal functor out of $([1],\mathrm{min},1)$. 

Further investing the equivalences
\[\mathrm{Alg}_{[1]^\mathrm{max}}(\mathcal O) \simeq \Ar(\mathrm{Alg}_{\Einf}(\mathcal O))\]
the non-unitally multiplicative and bimodule structures described above are in this perspective simply obtained via operad maps $\mathbb E_k^\mathrm{nu} \times [1]^\mathrm{max} \rightarrow \Ek \times [1]^\mathrm{min}$ and $\mathrm{BMod} \times [1]^\mathrm{max} \rightarrow \Ek \times [1]^\mathrm{min}$. The first one comes from the observation that the identity of $[1]$ is a non-unitally lax symmetric monoidal functor from $([1],\mathrm{max})$ to $([1],\mathrm{min})$. To construct the second map, note that it suffices to treat the case $k=1$, in which case all operads involved are ordinary. Denoting the colours of  $\mathrm{BMod}$ by $a_l,m,a_r$ and that of $\mathrm{Assoc} = \mathbb E_1$ by $a$, it is then easy to check that there is then a unique (planar, multi-coloured) operad map with 
\[\quad (m,0) \longmapsto (a,0), \quad (m,1) \longmapsto (a,1)  \quad \text{and} \quad  (a_l,0), (a_l,1), (a_r,0), (a_r,1) \longmapsto (a,1).\]

Of course, there are generalisation to operads other than $\Ek$; we leave the details to the reader.
\end{remark}

%E_\infty^nu ---> [1]^min

%{non-empty finite sets} partially defined maps with non-empty fibres.

%[1]^min = finite set S + f: S --> \{0,1\}       Hom(l,l') = \{\alpha mit for all i min_{k \in alpha^{-1}(i)} l_k \leq l'i}

%E_infty^nu x [1]^max--->   E_infty \times [1]^min
% (S,l)    ---- >  (S,l)

 In particular, if $A$ is a static ring, every ordinary (two-sided) ideal $I \subseteq A$ is a Smith-ideal in $A$ for $\C = \D(\mathbb Z)$, but note that for static $I$ the map $I \rightarrow A$ in a Smith-ideal does not need to be injective. The structure of a Smith ideal $I \rightarrow A$ is precisely what is required to equip the quotient map $A \rightarrow A/I$, with the structure of a ring map and in the stable setting this process actually yields an equivalence of categories.

\begin{proposition} Let $(\C,\otimes)$ be an $\mathbb E_k$-monoidal category. Assume that $\C$ admits all finite limits and colimits, $\C$ is pointed (i.e.~the initial and final object agree), and that $\otimes$ commutes with finite colimits in each variable.

For every $\mathbb E_k$-algebra $A$ in $\C$ extracting (co)fibres gives an adjunction
\[\cof \colon \mathrm{SmId}_{\Ek}(A) \leftrightarrows A/\mathrm{Alg}_{\mathbb E_k}(\C) \cocolon \fib.\]
If $\C$ is stable these are equivalences.
\end{proposition}

For $A$ a static ring we, in particular, find that ordinary ideals in $A$ are precisely those Smith-ideals $I$ in $A$ for which both $I$ and $A/I$ are static again. 

\begin{proof}
Consider $\Fun([1]^2,\C)$ equipped with Day convolution with respect to taking maxima in the first (say horizontal) component, and minima in the second (say vertical) component of $[1]^2$, and recall that Day convolution with respect to taking maxima in $[1]$ is just the pointwise monoidal structure on $\Ar(\C)$, and consider its full suboperad spanned by $\Fun_\ast([1]^2,\C)$, the subcategory of diagrams with vanishing lower right corner (i.e.\ the entry at $(1,0)$ is the initial and final object of $\C$), which is easily checked closed under the monoidal operation and contains the unit, and is thus $\Ek$-monoidal in its own right.

Now on the one hand, the further full suboperad thereof spanned by the cocartesian squares is equivalent to $\Fun([1]^\mathrm{min},\C^\otimes)^\mathrm{Day}$ by adding and removing the right column (i.e.\ the entries at $(1,0)$ and $(1,1)$): By functoriality of Day convolution the restriction map extends to one of operads, and one easily checks that the cocartesian squares are closed under Day convolution inside all squares (with lower right corner vanishing) and that the restriction is strong monoidal. Since the restriction functor is clearly an equivalence on underlying categories the claim follows. The inclusion of this subcategory into $\Fun_\ast([1]^2,\C)$ furthermore admits a right adjoint given by replacing the top right corner by the cofibre of the left vertical map. By \cite[Proposition 2.2.1.1]{HA}, this right adjoint then induces a lax $\Ek$-monoidal adjunction
\[\cof \colon \Fun([1]^\mathrm{min},\C^\otimes)^\mathrm{Day} \leftrightarrows \Fun_\ast([1]^\mathrm{max} \times [1]^\mathrm{min}, \C^\otimes)^\mathrm{Day} \cocolon \mathrm{res}_v,\]
adding and removing cofibres with strong $\Ek$-monoidal left adjoint.

On the other hand the subcategory of cartesian squares in $\Fun_\ast([1]^2,\C)$ is equivalent to $\Fun([1],\C)$ by forgetting the lower row (i.e.\ the entries at $(0,0)$ and $(1,0)$). This subcategory is not necessarily closed under the monoidal operation, but admits a left adjoint, given by replacing the lower left corner by the fibre of the top horizontal map. This left adjoint is easily checked to satisfy the conditions of \cite[Proposition 2.2.1.9]{HA}, which makes the subcategory strong $\Ek$-monoidal in its own right. In particular, it induces a lax $\Ek$-monoidal adjunction 
\[\mathrm{res}_h \colon \Fun_\ast([1]^\mathrm{max} \times [1]^\mathrm{min}, \C^\otimes)^\mathrm{Day} \leftrightarrows \Fun([1]^\mathrm{max},\C^\otimes)^\mathrm{Day} \cocolon \fib,\]
adding and removing fibres, again with strong $\Ek$-monoidal left adjoint.

Composing the two adjunctions (and taking fibres in the upper left corner) yields the first claim. If $\C$ is stable, then the subcategories of cartesian and cocartesian squares in $\Fun_\ast([1]^2,\C)$ agree, and restricting both adjunctions above to this subcategory gives the second claim.
\end{proof}

We end with two natural operations that can be performed on Smith ideals:

\begin{example}\label{smithex} There are Smith-ideal analogues of the sum and product of ordinary ideals. For this purpose let $\C$ be $\mathbb E_{k+l}$-monoidal with $l>0$ and let $I \rightarrow A$, $J \rightarrow B$ be two $\mathbb E_k$-Smith ideals in $\mathbb E_k$-algebras $A$ and $B$ in $\C$.
\begin{enumerate}
\item On the one hand, the category of $\mathbb E_k$-Smith ideals inherits a residual $\mathbb E_l$-monoidal structure and convolving two $\mathbb E_k$-Smith ideals provides the higher categorical way of taking the exterior sum $I \boxplus J \rightarrow A \otimes B$ of ideals, i.e.
\[(I \rightarrow A) \otimes^\mathrm{Day} (J \rightarrow B) \simeq (I \otimes B +_{I \otimes J} A \otimes J \rightarrow A \otimes B)\]
with corresponding $\mathbb E_k$-map $A \otimes B \longrightarrow A/I \otimes B/J$ and if $A = B$ convolving them over $0 \rightarrow A$ gives the interior sum $I + J \rightarrow A$, i.e.
\[(I \rightarrow A) \otimes^\mathrm{Day}_{(0 \rightarrow A)} (J \rightarrow A) \simeq (I +_{I \otimes_A J} J \rightarrow A)\]
with corresponding $\mathbb E_k$-map $A \rightarrow A/I \otimes_A A/J$.
\item On the other hand, taking pointwise tensor products generalises taking the product of ideals in that
\[(I \rightarrow A) \otimes^\mathrm{Day} (J \rightarrow B) \simeq (I \otimes J \rightarrow A \otimes B)\]
with corresponding $\mathbb E_k$-map $A \otimes B \rightarrow A/I \otimes B \times_{(A \otimes B)/(I \boxplus J)} A \otimes B/J$, with the interior version $I \cdot J \rightarrow A$ given by
\[(I \rightarrow A) \otimes^\mathrm{Day}_{(A \rightarrow A)} (J \rightarrow A) \simeq (I \otimes_A J \rightarrow A)\]
with corresponding map $A \rightarrow A/I \otimes_{A/I+J} A/J$; to see that this pointwise tensor product really preserves Smith-ideals, note that 
\[\Alg_{\mathbb E_{k}}\Fun^\mathrm{Day}([1]^{\mathrm{min}},-) \colon \Alg_{\Ek}(\Cat) \longrightarrow \Cat,\]
preserves products and so indeed lifts to a functor 
\[
\Alg_{\mathbb E_{k}}\Fun^\mathrm{Day}([1]^{\mathrm{min}},-) \colon \Alg_{\mathbb E_{k+l}}(\Cat) \longrightarrow \Alg_{\mathbb E_{l}}(\Cat)
\]
whose values are precisely equipped with the pointwise monoidal structures.
\end{enumerate}
\end{example}

\section{Proof of the main results}

For the proof recall the category $\Ek\text{-}\Mod(A)$ of $\Ek$-modules over $A$, as constructed in \cite[Sections 3.3 \& 3.4]{HA}, which is again $\Ek$-monoidal under $\otimes_A$. In contrast $\Mod(A)$ is only $\mathbb E_{k-1}$-monoidal; one has $\mathbb E_1\text{-}\Mod(A) = \mathrm{BiMod}(A,A)$ and $\mathbb E_\infty\text{-}\Mod(A) = \mathrm{Mod}(A)$.

\begin{proof}[Proof of Theorem B]
We start with the observation that $\mathrm{LQ}_A$ is indeed equivalent to a poset, i.e.\ its mapping spaces are either empty or contractible, by the characterisation of localisations of an $\mathbb E_k$-ring $A$ as $\otimes_A$-idempotent objects in $\mathbb E_k$-modules under $A$ (the analogue indeed holds for collections of idempotent objects in any monoidal category). Let us also immediately verify that $\ker(\pi_0\varphi)$ is indeed idempotent for $\varphi \colon A \rightarrow B$ a $\pi_0$-surjective localisation among connective $\mathbb E_k$-rings. Tensoring the fibre sequence $F \rightarrow A \rightarrow B$ with $F$ gives
\[
F \otimes_A F \longrightarrow F \longrightarrow F \otimes_A B
\]
and the right hand term vanishes since one has a fibre sequence
\[F \otimes_A B \longrightarrow A \otimes_A B \longrightarrow B \otimes_A B\]
whose right hand map (after identifying $A \otimes_A B \simeq B$) is a section of the multiplication $B \otimes_A B \rightarrow B$ and thus an equivalence. But $F$ is connective and the map $\pi_0(F) \rightarrow \ker\pi_0\varphi$ surjective by the long exact sequence of $\varphi$, whence a chase in the diagram
\[
\begin{tikzcd} \ker(\pi_0\varphi) \otimes_{\pi_0A} \ker(\pi_0\varphi) \ar[rr] && \ker(\pi_0\varphi) \\
\pi_0F \otimes_{\pi_0A} \pi_0F\ar[r,"\sim"] \ar[u] & \pi_0(F \otimes_A F) \ar[r,"\sim"] & \pi_0F \ar[u]
\end{tikzcd}
\]
shows that the multiplication $\ker(\pi_0\varphi) \otimes_{\pi_0A} \ker(\pi_0\varphi) \rightarrow \ker(\pi_0\varphi)$ is surjective as desired.

Next, we verify the last claim from the statement, i.e.\ that the inverse system $\widetilde I^{\otimes_A n}$ stabilises degreewise. In fact we show slightly more, namely that the cofibre $A/\widetilde I \otimes_A  \widetilde I^{\otimes_A n}$ of the canonical map $\widetilde I^{\otimes_A n+1} \rightarrow \widetilde I^{\otimes_A n}$ %induces an isomorphism on $\pi_i$ for $i<n-1$ and a surjection on $\pi_{n-1}$, or equivalently that its cofibre  
is $n$-connective.
%$A / \widetilde I^{\otimes_A n+1} \rightarrow A / \widetilde I^{\otimes_A n}$ is an isomorphism on $\pi_i$ for $i < n$ and surjective on $\pi_n$, or equivalently that its cofibre is $n$-connected. This will in particular imply that $A_I$ is connective with $\pi_0(A_I) = \pi_0(A)/I$.
%Computing the total cofibre of the square
%\[\begin{tikzcd}
%\widetilde I^{\otimes_A n+1} \ar[r] \ar[d] & \widetilde I^{\otimes_A n} \ar[d] \\
%A^{\otimes_A n+1} \ar[r,"\simeq"] & A^{\otimes_A n}
%\end{tikzcd}\]
%first horizontally then vertically and the other way around we find that this cofibre is $\mathbb S^1 \otimes A/\widetilde I \otimes_A  \widetilde I^{\otimes_A n}$. 
Since $A/\widetilde I = \GEM(\pi_0(A)/I)$ is an $\Ek$-ring annihilated by $I$, we immediately deduce that the homotopy groups of this cofibre are annihilated by $I$ (from both sides).

Now, for $n=0$, the connectivity claim is clear, and if we inductively assume that $A/\widetilde I \otimes_A  \widetilde I^{\otimes_A n}$ is $n$-connective, then 
\[
A/\widetilde I \otimes_A  \widetilde I^{\otimes_A n+1} = \left(A/\widetilde I \otimes_A  \widetilde I^{\otimes_A n}\right) \otimes_A \widetilde I
\]
is clearly also $n$-connective and its $n$th homotopy group is $\pi_n\left(A/\widetilde I \otimes_A  \widetilde I^{\otimes_A n}\right) \otimes_{\pi_0A} I$. Since the left hand term is annihilated by $I$, we compute
\[
\pi_n\left(A/\widetilde I \otimes_A  \widetilde I^{\otimes_A n}\right) \otimes_{\pi_0A} I= \pi_n\left(A/\widetilde I \otimes_A  \widetilde I^{\otimes_A n}\right) \otimes _{\pi_0(A)/I} \pi_0(A)/I \otimes_{\pi_0A} I
\]
\[
= \pi_n\left(A/\widetilde I \otimes_A  \widetilde I^{\otimes_A n}\right) \otimes_{\pi_0(A)/I} I/I^2 = 0.
\]

As the next step, we show that the tautological map $M = A \otimes_A M \rightarrow A/I^\infty \otimes_A M$ is an equivalence whenever the homotopy of $M$ is annihilated by $I$, or in other words that $I^\infty \otimes_A M \simeq 0$. We start with the simplest case $M= A/\widetilde I$, where the claim is equivalent to the multiplication map
\[
\left(\lim_{n\in \mathbb N^\op} \widetilde I^{\otimes_A n}\right) \otimes_A \widetilde I \longrightarrow \lim_{n\in \mathbb N^\op} \widetilde I^{\otimes_A n}
\]
being an equivalence. But since the limit stabilises degreewise and $\widetilde I$ is connective, we can move the limit out of the tensor product (the cofibre of the interchange map is a limit of terms with growing connectivity), and then the statement follows from finality.

For an arbitrary $A$-module $M$ concentrated in degree $0$ and killed by the action of $I$, choose a free resolution of $\pi_0M$ by $\pi_0(A)/I$-modules, which by the Dold-Kan theorem yields a diagram $F \colon \Delta^\op \rightarrow  \D(\pi_0(A)/I)$ with each $F_n$ concentrated in degree $0$, $\pi_0(F_n)$ free and $\colim_{\Delta^\op} F \simeq (\pi_0M)[0]$, so that $\colim_{\Delta^\op} \iota F \simeq M$, where $\iota$ is the composite $\D(\pi_0(A)/I) \simeq \Mod(A/\widetilde I) \rightarrow \Mod(A)$. But then 
\[
I^\infty \otimes_A M \simeq \colim_{k \in \Delta^\op} I^\infty \otimes_A \iota F_k \simeq 0
\]
since each $\iota F_k$ is a direct sum of $A/\widetilde I$. By exactness of $I^\infty \otimes_A (-)$, the claim then follows for each bounded $A$-module $M$ whose homotopy is annihilated by $I$ using the Postnikov tower of $M$. For bounded below $M$, we have
\[
I^\infty \otimes_A M \simeq I^\infty \otimes_A \left(\lim_{k \in \mathbb N^\op}\tau_{\leq k}M\right) \simeq \lim_{k \in \mathbb N^\op}I^\infty \otimes_A\tau_{\leq k}M \simeq 0
\]
by commuting the limit out using the same argument as above. Finally, for arbitrary $M$ whose homotopy is killed by $I$, we find
\[
I^\infty \otimes_A M \simeq I^\infty \otimes_A \left(\colim_{k \in \mathbb N}\tau_{\geq -k}M\right) \simeq \colim_{k \in \mathbb N}I^\infty \otimes_A\tau_{\geq -k}M \simeq 0.
\]

Now, since $A \rightarrow A/\widetilde I = \GEM(\pi_0(A)/I)$ is a map of $\Ek$-rings, it follows that $\widetilde I \rightarrow A$ is an $\Ek$-Smith-ideal in $A$. It then formally follows that so is $\widetilde I^{\otimes_A n} \rightarrow A$, whence $A/\widetilde I^{\otimes_A n}$ and thus $A/I^\infty$ are $\Ek$-rings by Example \ref{smithex}. Since $\pi_0(A/I^\infty) = \pi_0(A)/I$, all homotopy groups of $A/I^\infty$ are annihilated by $I$, and so the canonical map $A \rightarrow A/I^\infty$ induces an equivalence $A/I^\infty \rightarrow A/I^\infty \otimes_A A/I^\infty$, which shows that $A \rightarrow A/I^\infty$ is a localisation. Furthermore, it implies that the homotopy of every $A/I^\infty$-module is a $\pi_0(A)/I$-module, so combined with the previous point, we learn that the image of the fully faithful restriction functor $\Mod(A/I^\infty) \rightarrow \Mod(A)$ consists exactly of those modules whose homotopy is killed by $I$, as desired.

Finally, we are ready to verify that the construction $I \mapsto (A \rightarrow A/I^\infty)$ induces an inverse to taking kernels. The composition starting with an ideal is clearly the identity. So we are left to show that for every $\varphi \colon A \rightarrow B$ in $\mathrm{LQ}_A$ the canonical map $\psi \colon A/{\ker(\pi_0\varphi)}^\infty \rightarrow B$, arising from the homotopy of $B$ being annihilated by $\ker(\pi_0\varphi)$, is an equivalence. Per construction it induces an equivalence on $\pi_0$. By the lemma below, the functor $\psi_! = B \otimes_{A/\mathrm{ker}(\varphi)^\infty} - \colon \Mod(A/\ker(\pi_0\phi)^\infty) \rightarrow \Mod(B)$ is thus conservative when restricted to bounded below modules. But the map
\[
B \simeq \psi_!\left(A/\ker(\pi_0\varphi)^\infty\right) \xrightarrow{\psi_!(\varphi)} \psi_!(B) = B \otimes_{A/\ker(\pi_0\varphi)^\infty} B \simeq  B \otimes_A B
\]
is induced by the unit and thus an equivalence since $\varphi$ is a localisation. 
\end{proof}

\begin{lemma}
If $\psi \colon A \rightarrow B$ is a map of connective $\Eone$-rings which is an isomorphism on $\pi_0$, then
\[
B \otimes_A - \colon \Mod(A) \longrightarrow \Mod(B)
\]
is conservative when restricted to bounded below $A$-modules.
\end{lemma}

\begin{proof}
If $M \in \Mod(A)$ with $\pi_i(M) = 0$ for $i<n$, then $\pi_n(B \otimes_A M) = \pi_0(B) \otimes_{\pi_0(A)} \pi_n(M) = \pi_n(M)$, so if $M$ is bounded below with $B \otimes_A M \simeq 0$ then also $M \simeq 0$. Considering cofibres of morphisms, this implies the statement.
\end{proof}

Now, the main step in the deduction of Theorem A from Theorem B is to establish the structure of an animated commutative ring $R/I^\infty$ on $(\GEM R)/I^\infty$, whenever $R$ is itself an animated commutative ring. A simple argument for this is to note that the entire tower $(\GEM R)/\widetilde I^{\otimes_{\GEM R} -}$ inductively lifts to animated commutative rings on account of the cartesian squares
\[\begin{tikzcd}
(\GEM R)/\widetilde I^{\otimes_{\GEM R} n+1} \ar[r] \ar[d] & (\GEM R)/\widetilde I^{\otimes_{\GEM R} n} \ar[d] \\
(\GEM R)/\widetilde I \ar[r] & (\GEM R)/\widetilde I \otimes_{\GEM R}  (\GEM R)/\widetilde I^{\otimes_{\GEM R} n}
\end{tikzcd}\]
issuing from Example \ref{smithex} above, together with the equivalence $(\GEM R)/\widetilde I \simeq \pi_0(R)/I$. A slightly more conceptual way employs the Amitsur (or cobar) complex: 
%, in such a fashion that for a map of animated commutative rings $\psi \colon R \rightarrow S$ the induced map $\GEM R/\mathrm{ker}(\pi_0\varphi)^\infty \rightarrow \GEM S$ of $\mathbb E_\infty$-rings refines to a map of animated commutative rings. In this case, Theorem B provides the essential surjectivity of the functor out of $\mathrm{LQ}_R$ extracting the kernel on $\pi_0$ 
To this end, recall that for an algebra $B$ in a monoidal category $(\mathcal C,\otimes)$ it is the cosimplicial object in $\C$ with $[n] \mapsto B^{\otimes n+1}$ and face and degeneracy maps induced by the unit and multiplication, respectively, see e.g.\ \cite[Section 2.1]{MNN}. For a map $R \rightarrow S$ of animated commutative rings, we can consider it in $(R/\mathrm{AnCRing},\otimes^\mathbb L_R)$ and similarly for a map $A \rightarrow B$ of $\mathbb E_k$-rings we can consider it in $(A/\mathrm{Alg}_{{\mathbb E}_k}(\mathrm{Sp}),\otimes_A)$ and also in $(\mathbb E_k\text{-}\Mod(A),\otimes_A)$; these examples are connected by strong monoidal, limit preserving functors 
\[(R/\mathrm{AnCRing},\otimes^\mathbb L_R) \longrightarrow (\GEM R/\mathrm{Alg}_{{\mathbb E}_\infty}(\mathrm{Sp}),\otimes_{\GEM R}) \quad \text{and} \quad (A/\mathrm{Alg}_{{\mathbb E}_k}(\mathrm{Sp}),\otimes_A) \longrightarrow (\mathbb E_k\text{-}\Mod(A),\otimes_A).\]
In \cite[Proposition 2.14]{MNN} it is in particular shown, that the limit of the Amitsur complex for $A \rightarrow  \GEM (\pi_0(A)/I)$ formed in $(\mathbb E_k\text{-}\Mod(A),\otimes_A)$ agrees with $A/I^\infty$ as an $\mathbb E_k$-$A$-module. The characterisation of the $\mathbb E_k$-ring structure on $A/I^\infty$ as arising from being $\otimes_A$-idempotent then shows that this upgrades to an equivalence of $\mathbb E_k$-algebras under $A$. 

In particular, $(\GEM R)/I^\infty$ is the limit of the Amitsur complex of $\GEM R \rightarrow \GEM(\pi_0R/I)$ formed in $(\GEM R/\mathrm{Alg}_{{\mathbb E}_\infty}(\mathrm{Sp}),\otimes_{\GEM R})$, which allows us to lift this structure to that of an animated commutative ring by letting $R/I^\infty$ denote the limit of the Amitsur complex for $R \rightarrow \pi_0R/I$ in $(R/\mathrm{AnCRing},\otimes^\mathbb L_R)$. 

Regardless of the construction, we have $\GEM(R/I^\infty) \simeq (\GEM R)/I^\infty$ as $\mathbb E_\infty$-rings.

\begin{proof}[Proof of Theorem A]
We again observe that $\mathrm{LQ}_R$ is a poset, since $-\otimes^\mathbb L_R S \colon R/\mathrm{AnCRing} \rightarrow R/\mathrm{AnCRing}$ is a localisation onto its image for every $\varphi \colon R \rightarrow S$ in $\mathrm{LQ}_R$. The assignement $R \mapsto R/I^\infty$ thus gives a functor that is evidently right inverse to 
\[
\mathrm{LQ}_R \longrightarrow \{I \subseteq \pi_0R \mid I^2 = I\}, \quad \varphi \longmapsto \ker(\pi_0\varphi).
\]
Furthermore, from the case of $\mathbb E_\infty$-rings we learn that the natural map $S \simeq R \otimes_R^\mathbb L S \rightarrow R/I^\infty \otimes^\mathbb L_R S$ is an equivalence if and only if the homotopy groups of $S$ are annihilated by $I$. In this case we therefore obtain a map $R/I^\infty \rightarrow S$ of animated commutative rings under $R$, and in particular this applies in the case $I = \mathrm{ker}(\varphi)$. But by Theorem B there is only one map $\GEM R/\mathrm{ker}(\varphi)^\infty \rightarrow \GEM S$ under $\GEM R$ and this is an equivalence. Since the functor $\GEM \colon \mathrm{AnCRing} \rightarrow \mathrm{Alg}_{\mathbb E_\infty}(\mathrm{Sp})$ is conservative, we must thus also have $R/\mathrm{ker}(\varphi)^\infty \simeq S$ under $R$ as desired.
\end{proof}
\section{Examples and Remarks}

\begin{enumerate}
\item A different way of phrasing Theorem B is that there is a one-to-one correspondence between idempotent ideals in $\pi_0(A)$ and connective idempotent Smith-ideals in $A$ for every connective $\Ek$-algebra $A$, which makes $I \subseteq \pi_0(A)$ and $I^\infty \rightarrow A$ correspond. 
\item As a consequence of the classification of stable recollements, one obtains a cartesian square
\[\begin{tikzcd}
\mathrm{Mod}(A) \ar[rrrr,"{\hom_A(A/I^\infty,-) \Rightarrow A/I^\infty \otimes_A -}"] \ar[d] &&&&\mathrm{Ar}(\mathrm{Mod}(A/I^\infty)) \ar[d,"\mathrm{cof}"] \\
\mathrm{aMod}_I(A) \ar[rrrr,"{A/I^\infty \otimes_A \hom_A(I^\infty,-)}"] &&&& \mathrm{Mod}(A/I^\infty),
\end{tikzcd}\]
decomposing the module category of $A$ for every idempotent $I \subseteq \pi_0(A)$, and consequently exact squares
\[\begin{tikzcd}
M \ar[r] \ar[d] & \hom_A(I^\infty,M) \ar[d] && \hom_A(A/I^\infty,I^\infty \otimes_A M) \ar[r] \ar[d] & \hom(A/I^\infty,M) \ar[d] \\
M/I^\infty \ar[r] & \hom_A(I^\infty,M/I^\infty) &&  I^\infty \otimes_A M \ar[r] & M
\end{tikzcd}\]
for every $M \in \Mod(A)$, see \cite[Section A.5]{92}.
\item Either directly from the statement of the theorems, or via the construction using the Amitsur complex, one finds that for the exterior sum $I \boxplus_k J$ of two idempotent ideals $I \subseteq \pi_0(A)$ and $J \subseteq \pi_0(A')$ in two connective $k$-algebras $A$ and $A'$ ($k$ some $\mathbb E_2$-ring), that is the image of 
\[(\pi_0A \otimes_{\pi_0k} J) \oplus (I \otimes_{\pi_0 k} \pi_0A') \longrightarrow  \pi_0A \otimes_{\pi_0k} \pi_0A' = \pi_0(A \otimes^\mathbb L_k A'),\]
we have 
\[(A \otimes^{\mathbb L}_k A')/(I \boxplus_k J)^\infty \simeq A/I^\infty \otimes^{\mathbb L}_k A'/J^\infty,\]
or in other words $(I \boxplus_k J)^\infty$ is the exterior sum (over $k$) of the Smith-ideals $I^\infty$ in $A$ and $J^\infty$ in $A'$. 

This formula evidently also holds for three animated commutative rings in place of $k,A$ and $A'$. 
\item If $R$ is a static ring with an ideal $I$ that is flat as a left or right $R$-module and satisfies $I^2 = I$, then $I^{\otimes_{R}^\mathbb L n} = I^{\otimes_R n} = I$, so $I^\infty = I$ and $R/I^\infty = R/I$ is static.

As mentioned in the introduction, a commutative ring $R$ together with an idempotent, flat ideal $I \subseteq R$ is indeed one of the standard set-ups for almost mathematics, see e.g.\ \cite[Section 4]{Bhatt} for an exposition, and in this case $\mathrm{a}\mathcal D_I(R) \simeq \mathrm{aMod}_I(\GEM R)$ is the derived category of the ordinary category of almost $R$-modules. For example, this situation occurs whenever $(K,|\cdot|)$ is a perfectoid field: Then $\mathfrak m = \{x \in K \mid |x| <1\}$ is a flat and idempotent ideal in the valuation ring $\mathcal O = \{x \in K \mid |x|\leq 1\}$.
\item Let us also immediately note, that a finitely generated idempotent ideal $I$ in a (static) commutative ring $R$ is necessarily generated by single idempotent element $e$ by Nakayama's lemma and thus, as a direct summand, even projective over $R$. In this case $R/(e)^\infty = R/(e)$ is simply the factor of $R$ singled out by $e$, which can also be described as the ordinary localisation $R[(1-e)^{-1}]$ of $R$.
\item In fact, for $R$ static the animated ring $R/I^\infty$ is static if and only if $I \otimes_R^\mathbb L I \simeq I$ via the multiplication: The latter implies the former by the description of $\GEM R/I^\infty$ in Theorem B, and conversely if $R/I^\infty \simeq R/I$ we learn that $R/I \otimes_R^\mathbb L R/I \simeq R/I$, which by passing to fibres along the exact sequence $I \rightarrow R \rightarrow R/I$ first yields $I \otimes_R^\mathbb L R/I \simeq 0$ and then the claim. 

An example where this occurs without $I$ being flat is the ring $R_n = K[T_1^{1/2^\infty}, \dots, T_n^{1/2^\infty}]$ from the introduction with $I_n=(T_1^{1/2^\infty},\dots,T_n^{1/2^\infty})$. Then as a sequential colimit of principal ideals $I_1$ is flat over $R_1$, so the multiplicativity statement for exterior sums of ideals yields
\[R_n/I_n^\infty \simeq (R_1/I_1^\infty)^{\otimes^\mathbb L_K n} \simeq (R_1/I_1)^{\otimes^\mathbb L_K n} \simeq K^{\otimes_K^\mathbb L n} \simeq K.\]
But $I_n$ is no longer flat for $n \geq 2$: Setting $J_n = (T_1,\dots,T_n) \subseteq R_n$ we for example have
\[\mathrm{Tor}^{R_n}_i(I_n,R_n/J_n) \cong \begin{cases} I_n/J_n\cdot I_n & i = 0 \\ K^{n \choose {i+1}} & i \geq 1 \end{cases},\]
which can be read off from the exact sequence
\[I_n \otimes_{R_n}^\mathbb L R_n/J_n \longrightarrow R_n/J_n \longrightarrow K \otimes_{R_n}^\mathbb L R_n/J_n\]
in $\mathcal D(R_n/J_n)$ together with
\[K \otimes_{R_n}^\mathbb L R_n/J_n \simeq (K \otimes_{R_1}^\mathbb L R_1/T_1)^{\otimes^\mathbb L_K n} \simeq (\Lambda_K(K^{[1]}))^{\otimes^\mathbb L_K n} \simeq \Lambda_K((K^n)^{[1]}),\]
which in turn can be read off from the evident free resolution $R_1 \xrightarrow{T_1} R_1$ of $R_1/T_1$.

\item Whenever $I \otimes_R I$ is flat over $R$, one has $I^\infty \simeq I \otimes_R I$ (which does not generally agree with $I \otimes_R^\mathbb L I$ in this situation), and consequently $R/I^\infty \simeq R /^{\mathbb L} (I \otimes_R I)$ as mentioned in the introduction, where $/^{\mathbb L}$ denotes the cofibre in $\D(R)$, modelled by the commutative graded differential algebra with
\[(R /^{\mathbb L} (I \otimes_R I))_i = \begin{cases} R & i=0 \\ I \otimes_R I & i=1 \\ 0 & \text{else}\end{cases}\]
so that
\[
\pi_i(R/I^\infty)  = \begin{cases} R/I & i = 0 \\ \ker(I \otimes_R I \rightarrow I) & i = 1 \\ 0 & i \geq 2 \end{cases}
\]
in this case: The multiplication map $I^{\otimes_R^{\mathbb L} 2} \otimes^\mathbb L_R I^{\otimes_R^{\mathbb L} n} \rightarrow I^{\otimes_R^{\mathbb L} n}$
factors as
\[
I^{\otimes_R^{\mathbb L} 2} \otimes^\mathbb L_R I^{\otimes_R^{\mathbb L} n} \longrightarrow I^{\otimes_R 2} \otimes^\mathbb L_R I^{\otimes_R^{\mathbb L} n}  \longrightarrow  I^{\otimes_R^{\mathbb L} n}
\]
so the limit computing $I^\infty$ can be replaced by that over the terms $I^{\otimes_R 2} \otimes^\mathbb L_R I^{\otimes_R^{\mathbb L} n}$. But this system is constant, as can be seen inductively from the fibre sequence 
\[
 I^{\otimes_R 2} \otimes^\mathbb L_R I \longrightarrow I^{\otimes_R 2} \longrightarrow  I^{\otimes_R 2} \otimes^{\mathbb L}_R R/I 
\]
whose last term is the static module $I^{\otimes_R 2} \otimes_R  R/I$ by flatness of $I^{\otimes_R 2}$; and this equals $I \otimes_R I/I^2$ and thus vanishes as $I=I^2$.
\item Note that $I^{\otimes_R n} \cong I \otimes_R I$ for all $n \geq 2$ the moment $I$ is idempotent, e.g.\ by the stability assertion of Theorem B, so that no further flatness hypothesis can sensibly be put on tensor powers of $I$.
\item The condition $I\cdot\pi_nM = 0$ of $M \in \Mod(A)$ being almost zero is in fact equivalent to the a priori stronger condition that $I  \otimes_{\pi_0A} \pi_nM = 0$: For the former condition makes $\pi_nM$ into an $\pi_0(A)/I$-module so that
\[
I \otimes_{\pi_0A} \pi_n(M) = I \otimes_{\pi_0A} \pi_0(A)/I \otimes_{\pi_0(A)/I} \pi_nM = I/I^2 \otimes_{\pi_0(A)/I} \pi_nM = 0.
\]
%For example, this implies that for any idempotent ideal $I \subseteq R$, the multiplication map $I \otimes_R I \otimes_R I  \rightarrow I \otimes_R I$ is an isomorphism: One easily checks that the kernel $K$ of the multiplication $I \otimes_R I \rightarrow I$ satisfies $I \cdot K = 0$, since $i \otimes \sum_j r_j \otimes s_j = i \otimes (\sum_j r_js_j) = 0$ for all $i,r_j,s_j \in I$ such that $\sum_j r_j \otimes s_j \in K$. But then tensoring the short exact sequence 
%\[
%0 \longrightarrow K \longrightarrow I \otimes_R I \longrightarrow I \longrightarrow 0
%\]
%with $I$ yields an exact sequence
%\[
%I \otimes_R K \longrightarrow I \otimes_R I \otimes_R I \longrightarrow I \otimes_R I \longrightarrow 0,
%\]
%whose first term vanishes.   
%R(n) --- T^{2^{n+1}-2/2^{n+1}} ---> R(n) --- T^{1/2^{n}} ---> I(n)
%   |                                                          |
%incl                                              T^{1/2^{n+1}}
%   |                                                          |
%R(n+1) --- T^{2^{n+1}-1/2^{n+1} --> R(n+1) ---T^{1/2^{n+1}} ---> I(n+1)

%2^n+1-2/2^n+1 = X + 2^{n+1}-1/2^{n+1}   X = -1/2^n
\item In contrast to this, it need not be true, however, that $I \otimes_R^\mathbb L M \simeq 0$ for $M$ an $I$-almost zero $R$-module: For example, let $R = K[T^{1/2^\infty}]/T$ and $I = (T^{1/2^\infty})$, the ideal generated by all the $2$-power roots of $T$. Then $R/I = K$ is clearly almost $0$, but $I \otimes_{R}^\mathbb L K \simeq \bigoplus_{i\geq 1}K^{[2i-1]}$ does not vanish: Writing $ R(n) = K[T^{1/2^n}]/T$ and $I(n) = (T^{1/2^n})$ for the principal ideal therein, so that $R = \colim_n R(n)$ and $I = \colim_n I(n)$ and consequently $I \otimes^\mathbb L_R K = \colim_n I(n) \otimes^\mathbb L_{R(n)} K$, we can freely resolve the inclusion $I(n) \rightarrow I(n+1)$ by the periodic
\[
\begin{tikzcd}
\dots \ar[r] & R(n) \ar[rr,"\cdot T^{(2^n-1)/2^n}"] \ar[d,"\mathrm{incl}"] && R(n) \ar[rr,"\cdot T^{1/2^n}"] \ar[d,"\cdot T^{1/2^{n+1}}"] && R(n) \ar[rr,"\cdot T^{(2^n-1)/2^n}"] \ar[d,"\mathrm{incl}"] && R(n) \ar[d,"T^{1/2^{n+1}}"]\\
\dots \ar[r] & R(n+1) \ar[rr,"\cdot T^{(2^{n+1}-1)/2^{n+1}}"] && R(n+1) \ar[rr,"\cdot T^{1/2^{n+1}}"] && R(n+1) \ar[rr,"\cdot T^{(2^{n+1}-1)/2^{n+1}}"] && R(n+1).
\end{tikzcd}
\]
After tensoring with $K$ each term is $K$ with horizontal maps vanishing and vertical maps alternating between $0$ and $\id_K$. This gives the claim upon taking vertical colimits. 
\item The algebra $R = K[T^{1/2^\infty}]/T$ and ideal $I = (T^{1/2^\infty})/T$ from the previous point form a typical example for which $I \otimes_R I$, but not $I$ itself, is flat: We already used above that over $K[T^{1/2^\infty}]$ the ideal $(T^{1/2^\infty})$ is flat, and one easily checks that 
\[
I \otimes_{R} I \cong R \otimes_{K[T^{1/2^\infty}]} (T^{1/2^\infty}).
\]
%whereas $I$ itself is not the base change of $(T^{1/2^\infty})$. 
From the fibre sequence $I \otimes_R^\mathbb L I \rightarrow I \rightarrow I \otimes_R^\mathbb L K$ and the calculation in the previous point, one then reads off that $\ker(I \otimes_R I \rightarrow I) \cong K$, generated by $T^{1/2} \otimes T^{1/2}$, which gives
\[
\pi_\ast(R/I^\infty) = \Lambda_K(K^{[1]}),
\]
an exterior algebra on one generator in degree $1$, and also that there is an exact sequence
\[\bigoplus_{i \geq 1} K^{[2i]} \longrightarrow I \otimes_R^\mathbb L I \longrightarrow I^\infty\]
showing that the last two terms do not agree. 

Note also that in fact $R/I^\infty \simeq \Lambda_K(K^{[1]})$ as animated commutative $K$-algebras as the right hand side equals $\mathbb L\mathrm{Sym}_K^\ast(K^{[1]})$, the free animated commutative $K$-algebra on $K^{[1]}$, on account of the general formula $\mathbb L\mathrm{Sym}_A^n(M^{[1]}) \simeq \mathbb L\Lambda^n_A(M)^{[n]}$. They are not, however, equivalent as $R$-modules, let alone animated commutative $R$-algebras, since as an $R$-module the exterior algebra is just $K \oplus K^{[1]}$, which is not remotely idempotent over $R$ on account of the calculations from the previous item. 
\item For 
\[
R_n = K[T_1^{1/2^\infty}, \dots, T_n^{1/2^\infty}]/(T_1,\dots,T_n)
\]
and $I_n = (T_1^{1/2^\infty},\dots,T_n^{1/2^\infty})$, we then have $R_n/I_n^\infty = (R_1/I_1^\infty)^{\otimes_K^\mathbb L n}$ , so 
\[
\pi_\ast(R_n/I_n^\infty) = \Lambda_K(K^n[1]),
\]
by the previous point. This in particular shows that even for static rings $R$, the animated rings $R/I^\infty$ can have arbitrarily high non-trivial homotopy in the absence of any flatness assumption on $I$.

\item The categories $\mathrm{aMod}_I(A)$, alongside categories of sheaves on locally compact Hausdorff spaces, and categories of nuclear modules in condensed mathematics, are typical examples of compactly assembled categories that need not be compactly generated, e.g.\ for $I = (T^{1/2^\infty}) \subset K[T^{1/2^\infty}] = R$ this is due to Keller \cite{Keller}. In long anticipated work Efimov \cite{Efimov} recently defined a version of algebraic $K$-theory for such categories, and there results a fibre sequence
\[\mathrm K(\mathrm{aMod}_I(A)) \longrightarrow \mathrm K(A) \longrightarrow K(A/I^\infty)\]
of spectra for every idempotent $I \subseteq \pi_0A$. It is this connection to algebraic $K$-theory that originally sparked the present note.
\item For $R$ static and commutative we finally note that the category $\mathrm{a}\D(R,I)$ is compactly generated if and only if $I \subseteq R$ is a pure ideal, which by definition means that $R/I$ is flat over $R$ (which in turn implies that $I$ is flat and idempotent and so in particular $R/I \simeq R/I^\infty$).

If $I$ is pure then indeed $R/I \cong R[(1+I)^{-1}]$ as explained for example in \cite[Lemma 04PS]{stacks}, and for any $P \subseteq R$ the kernel of the extension by scalars functor $\D(R) \rightarrow \D(R[P^{-1}])$ is generated by the various $R/^{\mathbb L} p$ for $p \in P$, see e.g.\ \cite[Lemma 7.2.3.13]{HA}. For the converse let us say that a localisation $A \rightarrow S$ has perfectly generated fibre if the associated extension of scalars has compactly generated kernel. Then in \cite[Theorem 3.15]{thomason} Thomason showed (restricted to the affine case and slightly reinpreted) that
\[\mathrm{supp} \colon \{R \rightarrow S \mid \substack{\text{ derived localisation with }\\ \text{ perfectly generated fibre }}\} \leftrightarrows \{U \subseteq \mathrm{Spec}(R) \mid \substack{U \text{ is the intersection of its }\\ \text{ compact open neighbourhoods }}\}^\op \cocolon \mathbb R\Gamma(-,\mathcal O_R)\]
are inverse equivalences of posets, where 
\[\mathrm{supp}(S) = \{p \subseteq R \mid S \otimes^\mathbb L_R \kappa_p \simeq \kappa_p\}\]
with $\kappa_p$ the residue field at $p$, see also \cite{KP} for a simplified proof. 

In particular, any such $S$ is automatically coconnective and furthermore $S \otimes_R^\mathbb L - \colon \mathcal D(R) \rightarrow \mathcal D(S)$ preserves coconnective objects (by a simple induction on the representation of $\mathrm{supp}(S)$ through basic opens). If now $S$ is of the form $R/I^\infty$, it indeed follows that $S$ must be static and flat over $R$, hence $R/I^\infty \simeq R/I$ with $I$ pure.
\item Pure ideals $I$ that are not principal exist for example in every non-noetherian von Neumann ring (as then every ideal is pure). Explicit examples are $R = \prod_{j \in J} K_j$, $K_j$ fields and $I$ infinite with $I = \bigoplus_{j \in J} K_j$, in which case 
\[\bigoplus \colon \prod_{j \in J} \D(K_j) \longrightarrow \mathrm{a}\mathcal D(R,I)\]
is easily checked an equivalence, and the source is indeed generated by $\bigoplus_{j \in J} \Dperf(K_j)$. 

Another naturally occuring example is the augmention map $\mathbb Q[G] \rightarrow \mathbb Q$, where $G$ is an infinite torsion group; in this case $\mathbb Q[G] = \colim_{F \subseteq G} \mathbb Q[F]$, where $F$ exhausts the finite subgroups of $G$, displays $\mathbb Q[G]$ as a filtered colimit of semi-simple rings, which makes it a von Neumann ring. And if the augmentation ideal were principal, necessarily on an idempotent $e$, the kernel of multiplication with $e$ would map isomorphically onto $\mathbb Q$ via the augmentation. But $e$ necessarily lies in the image of some $\mathbb Q[F] \rightarrow \mathbb Q[A]$ so we have an injection
\[\mathrm{ker}(e \colon \mathbb Q[F] \rightarrow \mathbb Q[F]) \otimes_{\mathbb Q[F]} \mathbb Q[A] \longrightarrow \mathrm{ker}(e \colon \mathbb Q[A] \rightarrow \mathbb Q[A])\]
and the source is infinite dimensional if it does not vanish, and in the latter case $e=1$, which is absurd.

For an example outside the class of von Neumann rings, take the ideal $I$ of elements in $\mathrm{C}(J,\mathbb R) $ that vanish in a neighbourhood of some chosen $p \in J$, where $J \subseteq \mathbb R$ is an interval.

Note finally, that generally the augmentation map $R[G] \rightarrow R$, for $G$ a group, generally yields a derived localisation if and only if $\mathrm{H}_i^\mathrm{grp}(G,R)$ for all $i>0$, i.e.\ if $G$ is $R$-acyclic. In particular this yields lots of interesting non-commutative examples such as $\mathbb Z[V] \rightarrow \mathbb Z$ by \cite{szymikwahl}, where $V$ is the second Thompson group.  

\item Thomason's theorem also implies that the Frobenius endomorphism of an $\mathbb F_p$-algebra $R$ induces the identity on the poset of derived localisations with perfectly generated fibre, since this is certainly the case for the spectrum of $R$. Efimov recently posed the question whether something similar holds on the entire poset of derived localisations, that is on the smashing spectrum of $R$. The following is the counterexample already mentioned in the introduction: Let $S$ denote the set of finite strings of $0$'s and $1$'s and take
\[
R=\mathbb F_p[T_s\mid s \in S]/(T_s - T_{s \ast 0} \cdot T_{s \ast 1},T_s^p\mid s \in S)
\]
with $I$ generated by all the variables. Then under the Frobenius base change $- \otimes_R^\mathbb L \mathrm{Fr}^*R$ both $R \rightarrow R$ and $R \rightarrow R/I^\infty$ map to $R \rightarrow R$, and thus the map induced by Frobenius is not even injective: By Theorem A this can be tested on components, where we have to compare $R$ and $R/I \otimes_R \mathrm{Fr}^*R$, with the latter being the quotient of $R$ by the ideal generated by all the $p$-th powers of elements of $I$. But this is the trivial ideal.

\item Finally, Theorem A shows that an animated commutative ring structure on $R$ induces a unique compatible one on $R/I^\infty$; using the desciption of animated commutative rings as algebras over the monad of derived symmetric powers \cite[Section 4.2]{Raksit}, this implies that $\mathcal D(R)$ and $\mathcal D(R/I^\infty)$ are equipped with derived functors of $\mathrm{Sym}^n$, compatible under extension of scalars. It follows that also the category $\mathrm a\mathcal D_I(R)$ of derived almost modules carries such operations compatible with the left adjoint to the localisation $\mathcal D(R) \rightarrow \mathrm a\mathcal D_I(R)$. In most examples, these derived symmetric powers in fact simply descend from $\mathcal D(R)$ to $\mathrm a\mathcal D_I(R)$: This is true precisely if for every $k \geq 2$ the ideal $I$ is generated by the $k$-th powers of its elements and in particular it happens in all safe the previous examples above and always when $I\otimes_R I$ is flat; to see this combine \cite[Proposition 2.1.7 (ii)]{GR}, \cite[Theorem 14.1.57]{GR2}, \cite[Example 14.1.60]{GR2} and recall that $\mathbb L\mathrm{Sym}^n_R(M^{[2]}) \cong \mathbb L\Gamma^n_R(M)^{[2n]}$, so that ruling out descent for derived symmetric powers is the same as ruling it out for derived divided powers. The ring from the previous example is an explicit case where the statement fails.

This is the only important structural result we are aware of, that actually requires a flatness assumption.
\end{enumerate}

\bibliographystyle{amsalpha} 
 
\newcommand{\etalchar}[1]{$^{#1}$}
\providecommand{\bysame}{\leavevmode\hbox to3em{\hrulefill}\thinspace}
\providecommand{\MR}{\relax\ifhmode\unskip\space\fi MR }
\providecommand{\MRhref}[2]{%
  \href{http://www.ams.org/mathscinet-getitem?mr=#1}{#2}
}
\providecommand{\href}[2]{#2}

\end{document}